\documentclass[12pt]{article}
\usepackage{amssymb,amsmath}
\usepackage{hyperref}
\usepackage{graphicx}
\usepackage{color}

\begin{document}

\title{\LARGE\bf An alternative representation of the Vi\'{e}te\text{'}s formula for pi by Chebyshev polynomials of the first kind}

\author{
\normalsize\bf S. M. Abrarov\footnote{\scriptsize{Dept. Earth and Space Science and Engineering, York University, Toronto, Canada, M3J 1P3.}}\, and B. M. Quine$^{*}$\footnote{\scriptsize{Dept. Physics and Astronomy, York University, Toronto, Canada, M3J 1P3.}}}

\date{September 19, 2016}
\maketitle

\begin{abstract}
There are several reformulations of the Vi\'{e}te\text{'}s formula for pi that have been reported in the modern literature. In this paper we show another analog to the Vi\'{e}te\text{'}s formula for pi by Chebyshev polynomials of the first kind.
\vspace{0.25cm}
\vspace{0.25cm}
\end{abstract}

\section{Introduction}
The sinc function, also known as the cardinal sine function, is defined as \cite{Gearhart1990, Kac1959}
\[
\text{sinc}\left( t \right)=\left\{ \begin{aligned}
\frac{\sin \left( t \right)}{t}, &\qquad\qquad t \ne 0 \\ 
1, &\qquad\qquad t = 0. \\ 
\end{aligned} 
\right.
\]
The sinc function finds many applications in sampling, spectral methods, differential equations and numerical integration \cite{Gearhart1990, Stenger2011, Rybicki1989, Lether1998, Quine2013, Abrarov2015, Ortiz-Gracia2016}.

More than four centuries ago the French lawyer and amateur mathematician Fran\c{c}ois Vi\'{e}te found a fabulous relation showing how the sinc function can be represented elegantly as an infinite product of the cosines \cite{Gearhart1990, Kac1959}
\begin{equation}\label{eq_1}
\text{sinc}\left( t \right)=\prod\limits_{m=1}^{\infty }{\cos \left( \frac{t}{{{2}^{m}}} \right)}.
\end{equation}
Since
$$
\text{sinc}\left( \frac{\pi }{2} \right)=\frac{2}{\pi }
$$
we may attempt to substitute the argument $t=\pi /2$ into right side of equation \eqref{eq_1}. Thus, using repeatedly for each $m$ the following cosine identity for double angle
$$
\cos \left( 2\theta_m  \right)=2{{\cos }^{2}}\left( \theta_m  \right)-1
$$
or
$$
\cos \left( \theta_m  \right)=2{{\cos }^{2}}\left( \theta_{m+1}  \right)-1 \Leftrightarrow \cos \left( {{\theta }_{m+1}} \right)=\sqrt{\frac{\cos \left( {{\theta }_{m}} \right)+1}{2}},
$$
where
$$
\theta_m =\frac{\pi /2}{{{2}^{m}}}, \quad \theta_{m+1} = \frac{\theta_{m}}{2}
$$
and taking into account that
$$
\cos \left( \theta_{1} \right) = \cos \left( \frac{\pi /2}{{{2}^{1}}} \right)=\frac{\sqrt{2}}{2},
$$
we can find the following sequence
$$
\cos \left( \frac{\pi /2}{{{2}^{2}}} \right)=\frac{\sqrt{2+\sqrt{2}}}{2},
$$
$$
\cos \left( \frac{\pi /2}{{{2}^{3}}} \right)=\frac{\sqrt{2+\sqrt{2+\sqrt{2}}}}{2},
$$
$$
\vdots
$$
\begin{equation}\label{eq_2}
\cos \left( \frac{\pi /2}{{{2}^{m}}} \right)=\frac{\overbrace{\sqrt{2+\sqrt{2+\sqrt{2+\cdots +\sqrt{2}}}}}^{m\,\,\text{square}\,\,\text{roots}}}{2}.
\end{equation}
Therefore, from the equations \eqref{eq_1} and \eqref{eq_2} we obtain the Vi\'{e}te\text{'}s infinite product formula for the constant pi (in radicals consisting of square roots and twos only) \cite{Osler1999, Servi2003, Levin2005, Levin2006, Kreminski2008}
\[
\begin{aligned}
\text{sinc}\left( \frac{\pi}{2} \right) &=\cos \left( \frac{\pi /2}{{{2}^{1}}} \right)\cos \left( \frac{\pi /2}{{{2}^{2}}} \right)\cos \left( \frac{\pi /2}{{{2}^{3}}} \right)\cdots  \\ 
 & =\frac{\sqrt{2}}{2}\frac{\sqrt{2+\sqrt{2}}}{2}\frac{\sqrt{2+\sqrt{2+\sqrt{2}}}}{2}\cdots = \frac{2}{\pi }
\end{aligned}
\]
that can be conveniently rewritten as
\begin{equation}\label{eq_3}
\frac{2}{\pi }=\underset{M\to \infty }{\mathop{\lim }}\,\prod\limits_{m=1}^{M}{\frac{{{a}_{m}}}{2}},
\end{equation}
where ${{a}_{m}}=\sqrt{2+{{a}_{m-1}}}$ and ${{a}_{1}}=\sqrt{2}$.

Several reformulations of the Vi\'{e}te\text{'}s formula \eqref{eq_3} for pi have been reported in the modern literature \cite{Osler1999, Servi2003, Levin2005, Levin2006, Kreminski2008}. Notably, Osler has shown by \text{``}double product\text{''} generalization a direct relationship between the classical Vi\'{e}te\text{'}s and Wallis\text{'}s infinite products for pi (see equation (3) in \cite{Osler1999}). In this work we derive another equivalent to the Vi\'{e}te\text{'}s formula for pi expressed in terms of the Chebyshev polynomials of the first kind.

\section{Derivation}

The Chebyshev polynomials ${{T}_{m}}\left( x \right)$ of the first kind can be defined by the following recurrence relation \cite{Press1992, Mathews1999, Zwillinger2012}
$$
{{T}_{m}}\left( x \right)=2x{{T}_{m-1}}\left( x \right)-{{T}_{m-2}}\left( x \right),
$$
where ${{T}_{1}}\left( x \right)=x$ and ${{T}_{0}}\left( x \right)=1$. It should be noted that the recurrence procedure is not required in computation since these polynomials can also be determined directly by using, for example, a simple identity
\[
T_{m}\left( x \right)=x^{m}\sum\limits_{n=0}^{\left\lfloor m/2 \right\rfloor } \binom{m}{2n}{{\left( {1 - {x}^{-2}} \right)}^{n}}.
\]

Due to a remarkable property of the Chebyshev polynomials
$$
\cos \left( m \alpha \right) = {{T}_{m}}\left( \cos \left( \alpha \right) \right),
$$
making change of the variable in form
$$
\alpha = \cos \left( \frac{t}{{{2}^{M}}} \right)
$$
results in
\begin{equation}\label{eq_4}
\cos \left( \frac{2m-1}{{{2}^{M}}}t \right)={{T}_{2m-1}}\left( \cos \left( \frac{t}{{{2}^{M}}} \right) \right).
\end{equation}
Consequently, substituting equation \eqref{eq_4} into the following product-to-sum identity \cite{Quine2013, Abrarov2015, Ortiz-Gracia2016}
\begin{equation}\label{eq_5}
\prod\limits_{m=1}^{M}{\cos \left( \frac{t}{{{2}^{m}}} \right)}=\frac{1}{{{2}^{M-1}}}\sum\limits_{m=1}^{{{2}^{M-1}}}{\cos \left( \frac{2m-1}{{{2}^{M}}}t \right)}
\end{equation}
yields
\begin{equation}\label{eq_6}
\prod\limits_{m=1}^{M}{\cos \left( \frac{t}{{{2}^{m}}} \right)}=\frac{1}{{{2}^{M-1}}}\sum\limits_{m=1}^{{{2}^{M-1}}}{{{T}_{2m-1}}\left( \cos \left( \frac{t}{{{2}^{M}}} \right) \right)}.
\end{equation}
It can be shown that the right side of the equation \eqref{eq_6} can be further simplified and represented by a single Chebyshev polynomial of the second kind (see {\it{Appendix A}}).

Comparing equations \eqref{eq_1} and \eqref{eq_5} we can see that the infinite product of cosines for the sinc function can be transformed into infinite sum of cosines \cite{Abrarov2015}
\begin{equation}\label{eq_7}
\text{sinc}\left( t \right)=\underset{M\to \infty }{\mathop{\lim }}\,\frac{1}{{{2}^{M-1}}}\sum\limits_{m=1}^{{{2}^{M-1}}}{\cos \left( \frac{2m-1}{{{2}^{M}}}t \right)}.
\end{equation}
Since the right side of equation \eqref{eq_5} represents a truncation of the limit \eqref{eq_7} by a finite value of upper integer ${{2}^{M-1}}$ in summation, it is simply the incomplete cosine expansion of the sinc function. Indeed, if the condition ${{2}^{M-1}}>>1$ is satisfied, then the incomplete cosine expansion of the sinc function quite accurately approximates the original sinc function as given by \cite{Abrarov2015, Ortiz-Gracia2016}
$$
\frac{1}{{{2}^{M-1}}}\sum\limits_{m=1}^{{{2}^{M-1}}}{\cos \left( \frac{2m-1}{{{2}^{M}}}t \right)}\approx \text{sinc}\left( t \right).
$$

It is interesting to note that comparing equations \eqref{eq_1} and \eqref{eq_6} we can also write now
\footnotesize
$$
\text{sinc}\left( t \right)=\underset{M\to \infty }{\mathop{\lim }}\,\frac{1}{{{2}^{M-1}}}\sum\limits_{m=1}^{{{2}^{M-1}}}{{{T}_{2m-1}}\cos \left( \frac{t}{{{2}^{M}}} \right)}=\underset{M\to \infty }{\mathop{\lim }}\,\frac{1}{{{2}^{M-1}}}\sum\limits_{m=1}^{{{2}^{M-1}}}{{{T}_{2m-1}}\cos \left( \frac{t/2}{{{2}^{M-1}}} \right)}
$$
\normalsize
or
$$
\text{sinc}\left( t \right)=\underset{L\to \infty }{\mathop{\lim }}\,\frac{1}{L}\sum\limits_{\ell =1}^{L}{{{T}_{2\ell -1}}\cos \left( \frac{t}{2L} \right)},
$$
since we can imply that $L={{2}^{M-1}}$.

At $t=\pi /2$ the equation \eqref{eq_6} provides
\begin{equation}\label{eq_8}
\prod\limits_{m=1}^{M}{\cos \left( \frac{\pi /2}{{{2}^{m}}} \right)}=\frac{1}{{{2}^{M-1}}}\sum\limits_{m=1}^{{{2}^{M-1}}}{{{T}_{2m-1}}\left( \cos \left( \frac{\pi /2}{{{2}^{M}}} \right) \right)}.
\end{equation}
Applying equation \eqref{eq_2} again for each $m$ repeatedly, the product-to-sum identity \eqref{eq_8} can be rearranged in form
\footnotesize
\[
\begin{aligned}
\frac{\sqrt{2}}{2}\frac{\sqrt{2+\sqrt{2}}}{2}\frac{\sqrt{2+\sqrt{2+\sqrt{2}}}}{2} & \ldots \frac{\overbrace{\sqrt{2+\sqrt{2+\sqrt{2+\cdots +\sqrt{2}}}}}^{M\,\,\text{square}\,\,\text{roots}}}{2}\\
& = \frac{1}{{{2}^{M-1}}}\sum\limits_{m=1}^{{{2}^{M-1}}}{{{T}_{2m-1}}\left( \frac{\overbrace{\sqrt{2+\sqrt{2+\sqrt{2+\cdots +\sqrt{2}}}}}^{M\,\,\text{square}\,\,\text{roots}}}{2} \right)}
\end{aligned}
\]
\normalsize
or
\begin{equation}\label{eq_9}
\prod\limits_{m=1}^{M}{\frac{{{a}_{m}}}{2}}=\frac{1}{{{2}^{M-1}}}\sum\limits_{m=1}^{{{2}^{M-1}}}{{{T}_{2m-1}}\left( \frac{{a}_{M}}{2} \right)}.
\end{equation}

Increase of the integer $M$ approximates the product of cosines on the left side of equation \eqref{eq_9} closer to the value $2/\pi$. This signifies that the right side of the equation \eqref{eq_9} also tends to $2/\pi $ as the integer $M$ increases. Consequently, this leads to
$$
\frac{2}{\pi }=\underset{M\to \infty }{\mathop{\lim }}\,\frac{1}{{{2}^{M-1}}}\sum\limits_{m=1}^{{{2}^{M-1}}}{{{T}_{2m-1}}\left( \frac{\overbrace{\sqrt{2+\sqrt{2+\sqrt{2+\cdots +\sqrt{2}}}}}^{M\,\,\text{square}\,\,\text{roots}}}{2} \right)}
$$
or
\begin{equation}\label{eq_10}
\frac{2}{\pi }=\underset{M\to \infty }{\mathop{\lim }}\,\frac{1}{{{2}^{M-1}}}\sum\limits_{m=1}^{{{2}^{M-1}}}{{{T}_{2m-1}}\left( \frac{{a}_{M}}{2} \right)}.
\end{equation}

The equation \eqref{eq_10} is completely identical to the Vi\'{e}te infinite product \eqref{eq_3} for the constant pi. Since the relation \eqref{eq_8} remains valid for any integer $M$, the equation \eqref{eq_10} can be regarded as a product-to-sum transformation of the Vi\'{e}te\text{'}s formula for pi. 

It should be noted that the equation \eqref{eq_10} can be readily rearranged as a single Chebyshev polynomial of the second kind (see {\it{Appendix B}}).

\section{Conclusion}
We show a new analog to the Vi\'{e}te\text{'}s formula for pi represented in terms of the Chebyshev polynomials of the first kind. This approach is based on a product-to-sum transformation of the Vi\'{e}te\text{'}s formula.

\section*{Acknowledgments}

This work is supported by National Research Council Canada, Thoth Technology Inc. and York University.

\section*{Appendix A}

The Chebyshev polynomials ${{U}_{m}}\left( x \right)$ of the second kind can also be defined by the recurrence relation. Specifically, we can write \cite{Zwillinger2012}
$$
{{U}_{m}}\left( x \right)=2x{{U}_{m-1}}\left( x \right)-{{U}_{m-2}}\left( x \right),
$$
where ${{U}_{1}}\left( x \right)=2x$ and ${{U}_{0}}\left( x \right)=1$.

There is a simple relation for sum of the odd Chebyshev polynomials of the first kind
\[
\tag{A.1}\label{A.1}
{{U}_{K}}\left( x \right)=2\sum\limits_{k\,\,\text{odd}}^{K}{{{T}_{k}}\left( x \right)},
\]
where $K$ is an odd integer. Consequently, substituting equation \eqref{eq_6} into relation \eqref{A.1} provides\footnote{The subscript $2^M-1$ should not be confused with notation $2^{M-1}$ that has been used in some equations earlier.}
\[
\prod\limits_{m=1}^{M}{\cos \left( \frac{t}{{{2}^{m}}} \right)}=\frac{1}{{{2}^{M}}}{{U}_{{{2}^{M}}-1}}\left( \cos \left( \frac{t}{{{2}^{M}}} \right) \right).
\]
According to equation \eqref{eq_1} tending $M$ to infinity leads to the limit
\[
\text{sinc}\left( t \right)=\underset{M\to \infty }{\mathop{\lim }}\,\frac{1}{{{2}^{M}}}{{U}_{{{2}^{M}}-1}}\left( \cos \left( \frac{t}{{{2}^{M}}} \right) \right)
\]
or
\[
\tag{A.2}\label{A.2}
\text{sinc}\left( t \right)=\underset{N\to \infty }{\mathop{\lim }}\,\frac{1}{N}{{U}_{N-1}}\left( \cos \left( \frac{t}{N} \right) \right),  
\]
since we can imply that $2^M = N$. Obviously at $N >> 1$, one can truncate equation \eqref{A.2} to approximate the sinc function by a single Chebyshev polynomial of the second kind as
\[
\text{sinc}\left( t \right) = \frac{1}{N}{{U}_{N-1}}\left( \cos \left( \frac{t}{N} \right) \right)+\epsilon \left( t \right),
\]
where $\epsilon \left( t \right)$ is the error term. For example, taking $N=16$ results in
\[
\begin{aligned}
\text{sinc}\left(t\right) = & \,\, 2048 \cos ^{15}\left(\frac{t}{16}\right)-7168 \cos ^{13}\left(\frac{t}{16}\right)+9984 \cos ^{11}\left(\frac{t}{16}\right)\\
&-7040 \cos ^9\left(\frac{t}{16}\right)+2640 \cos ^7\left(\frac{t}{16}\right)-504 \cos ^5\left(\frac{t}{16}\right)\\
&+42 \cos ^3\left(\frac{t}{16}\right)-\cos \left(\frac{t}{16}\right) + \epsilon \left( t \right),
\end{aligned}
\]
where within the range $-10 \leq t \leq 10$ the error term satisfies $\left| \epsilon \left( t \right) \right| < 0.006$. As we can see, this approach quite accurately approximates the sinc function even if the integer $N$ in the limit \eqref{A.2} is not very large.

\section*{Appendix B}

Substituting equation \eqref{eq_10} into relation \eqref{A.1} we can express the Vi\'{e}te\text{'}s formula for pi by a single Chebyshev polynomial of the second kind as given by
\[
\frac{2}{\pi}=\underset{M\to \infty }{\mathop{\lim }}\,\frac{1}{{{2}^{M}}}{{U}_{{{2}^{M}}-1}}\left( \frac{\overbrace{\sqrt{2+\sqrt{2+\sqrt{2+\cdots \sqrt{2}}}}}^{M\,\,\text{square}\,\,\text{roots}}}{2} \right)
\]
or
\[
\tag{B.1}\label{B.1}
\frac{2}{\pi}=\underset{M\to \infty }{\mathop{\lim }}\,\frac{1}{{{2}^{M}}}{{U}_{{{2}^{M}}-1}}\left( \frac{{{a}_{M}}}{2} \right).
\]

Although the equation \eqref{B.1} is more simple, the equation \eqref{eq_10} reflects explicitly the product-to-sum transformation of the Vi\'{e}te\text{'}s formula \eqref{eq_3} for the constant $\pi$.

\bigskip

\end{document}